\documentclass[a4paper,12pt]{amsart}
\usepackage{amsfonts,amsmath,amsthm,amssymb}
\usepackage{enumitem}
\usepackage{mathrsfs}
\usepackage[margin=1in]{geometry}
\usepackage[colorlinks=true,linkcolor=blue,filecolor=blue,citecolor=red]{hyperref}

\usepackage{cleveref}

\theoremstyle{plain}
\newtheorem{theorem}{Theorem}[section]

\theoremstyle{definition}

\begin{document}
\title[Congruences for generalized cubic partitions modulo primes]{A note on congruences for generalized cubic partitions modulo primes}
\author[Russelle Guadalupe]{Russelle Guadalupe}
\address{Institute of Mathematics, University of the Philippines, Diliman\\
	Quezon City 1101, Philippines}
\email{rguadalupe@math.upd.edu.ph}

\renewcommand{\thefootnote}{}



\renewcommand{\thefootnote}{\arabic{footnote}}
\setcounter{footnote}{0}

\begin{abstract}
Recently, Amdeberhan, Sellers, and Singh introduced the notion of a generalized cubic partition function $a_c(n)$ and proved two isolated congruences via modular forms, namely, $a_3(7n+4)\equiv 0\pmod{7}$ and $a_5(11n+10)\equiv 0\pmod{11}$. In this paper, we provide another proof of these congruences by using classical $q$-series manipulations. We also give infinite families of congruences for $a_c(n)$ for primes $p\not\equiv 1\pmod{8}$.
\end{abstract}

\maketitle

\section{Introduction}
Throughout this paper, we denote $f_m = \prod_{n\geq 1}(1-q^{mn})$ for a positive integer $m$ and a complex number $q$ with $|q| < 1$. Recall that a partition of a positive integer $n$ is a nonincreasing sequence of positive integers, known as its parts, whose sum is $n$. A cubic partition of $n$ is a partition of $n$ whose even parts may appear in two colors. We denote $a(n)$ as the number of cubic partitions of $n$ and set $a(0):=1$. Then the generating function of $a(n)$ is 
\[\sum_{n=0}^\infty a(n)q^n=\dfrac{1}{f_1f_2}.\]
Chan \cite{chan} showed that $a(3n+2)\equiv 0\pmod{3}$ by establishing the remarkable identity
\[\sum_{n=0}^\infty a(3n+2)q^n = 3\dfrac{f_3^3f_6^3}{f_1^4f_2^4},\]
which is an analogue of the identity of Ramanujan \cite[pp. 210--213]{ramp} given by
\[\sum_{n=0}^\infty p(5n+4)q^n = 5\dfrac{f_5^5}{f_1^6},\]
where $p(n)$ is number of partitions of $n$.

Recently, Amdeberhan, Sellers, and Singh \cite{amdb} introduced the notion of a generalized cubic partition of $n$, which is a partition of $n$ whose even parts may appear in $c\geq 1$ different colors. We denote $a_c(n)$ as the number of such generalized cubic partitions of $n$ and set $a_c(0):=1$. Then the generating function of $a_c(n)$ is 
\[\sum_{n=0}^\infty a_c(n)q^n=\dfrac{1}{f_1f_2^{c-1}}.\]
Using the theory of modular forms, they proved the following congruences for $a_c(n)$ modulo $7$ and $11$.
\begin{theorem}[\cite{amdb}]\label{thm11}
	For all $n\geq 0$,
	\begin{align}
		a_3(7n+4)&\equiv 0\pmod{7},\label{eq11}\\
		a_5(11n+10)&\equiv 0\pmod{11}.\label{eq12}
	\end{align}
\end{theorem}
In this paper, we offer another proof of Theorem \ref{thm11} using classical $q$-series manipulations. We also provide infinite families of congruences for $a_c(n)$ modulo primes $p\not\equiv 1\pmod{8}$.

The rest of the paper is organized as follows. In Section \ref{sec2}, we present another proof of Theorem \ref{thm11} using the identities of Euler and Ramanujan. In Section \ref{sec3}, we prove two infinite families of congruences for $a_c(n)$ modulo primes $p\not\equiv 1\pmod{8}$ using another identity of Ramanujan and the result of Ahlgren \cite{ahl}.

\section{Another proof of Theorem \ref{thm11}}\label{sec2}

\begin{proof}
	To prove Theorem \ref{thm11}, we recall the identity of Euler \cite[(1.7.1)]{hirsc}
	\begin{equation}\label{eq21}
		f_1=\sum_{n=-\infty}^\infty (-1)^nq^{n(3n+1)/2}
	\end{equation}
	and the identity of Ramanujan \cite[(10.7.3)]{hirsc}
	\begin{equation}\label{eq22}
		\dfrac{f_1^5}{f_2^2}=\sum_{n=-\infty}^\infty (-1)^n(6n+1)q^{n(3n+1)/2}.
	\end{equation}
	
	By the binomial theorem, we have $f_7\equiv f_1^7\pmod{7}$, so from equations (\ref{eq21}) and (\ref{eq22}), 
	
	\begin{align}
		\sum_{n=0}^\infty a_3(n)q^n&=\dfrac{1}{f_1f_2^2}\equiv \dfrac{1}{f_7}\cdot\dfrac{f_1^5}{f_2^2}\cdot f_1\nonumber\\
		&\equiv \dfrac{1}{f_7}\sum_{m,\,n=-\infty}^\infty (-1)^{m+n}(6m+1)q^{(3m^2+m)/2+(3n^2+n)/2}\pmod{7}.\label{eq23}
	\end{align}
	We consider the equation 
	\[\dfrac{3m^2+m}{2}+\dfrac{3n^2+n}{2}\equiv 4\pmod{7},\]
	which is equivalent to
	\begin{align}
		(6m+1)^2+(6n+1)^2\equiv 0\pmod{7}.\label{eq24}
	\end{align}
	Since $7\equiv 3\pmod{4}$, $-1$ is a quadratic nonresidue modulo $7$, so the solution of equation (\ref{eq24}) is $6m+1\equiv6n+1\equiv 0\pmod{7}$. Thus, extracting the terms containing $q^{7n+4}$ on both sides of equation (\ref{eq23}), dividing by $q^4$, and then replacing $q^7$ with $q$, we get equation (\ref{eq11}).
	
	On the other hand, we have $f_{11}\equiv f_1^{11}\pmod{11}$, so equation (\ref{eq22}) implies that
	\begin{align}
		\sum_{n=0}^\infty a_5(n)q^n&=\dfrac{1}{f_1f_2^4}\equiv \dfrac{1}{f_{11}}\cdot\left(\dfrac{f_1^5}{f_2^2}\right)^2\nonumber\\
		&\equiv \dfrac{1}{f_{11}}\sum_{m,\,n=-\infty}^\infty (-1)^{m+n}(6m+1)(6n+1)q^{(3m^2+m)/2+(3n^2+n)/2}\pmod{11}.\label{eq25}
	\end{align}
	
	We consider the equation 
	\[\dfrac{3m^2+m}{2}+\dfrac{3n^2+n}{2}\equiv 10\pmod{11},\]
	which is equivalent to
	\begin{align}
		(6m+1)^2+(6n+1)^2\equiv 0\pmod{11}.\label{eq26}
	\end{align}
	Since $11\equiv 3\pmod{4}$, $-1$ is a quadratic nonresidue modulo $11$, so the solution of equation (\ref{eq26}) is $6m+1\equiv6n+1\equiv 0\pmod{11}$. Thus, extracting the terms containing $q^{11n+10}$ on both sides of equation (\ref{eq25}), dividing by $q^{10}$, and then replacing $q^{11}$ with $q$, we arrive at equation (\ref{eq12}). 
\end{proof}

\section{Congruences for $a_c(n)$ modulo primes $p\not\equiv 1\pmod{8}$}\label{sec3}

We now prove two infinite families of congruences for $a_c(n)$ modulo primes $p\not\equiv 1\pmod{8}$. We first give the following result for primes $p\equiv 5,7\pmod{8}$, which generalizes equation (\ref{eq11}) in Theorem \ref{thm11}.

\begin{theorem}\label{thm31}
	Let $p\equiv 5, 7\pmod{8}$ be a prime and $0\leq l\leq p-1$ be a nonnegative integer with $p\mid 8l+3$. Then for all $n\geq 0$, 
	\begin{equation}\label{eq31}
		a_{p-4}(pn+l)\equiv 0\pmod{p}.
	\end{equation}
\end{theorem}

\begin{proof}
	We start with the following identity of Ramanujan \cite[(10.7.7)]{hirsc}
	\begin{equation}\label{eq32}
		\dfrac{f_2^5}{f_1^2}=\sum_{n=-\infty}^\infty (-1)^n(3n+1)q^{3n^2+2n}.
	\end{equation}
	With $f_{2p}\equiv f_2^p\pmod{p}$, we see from equations (\ref{eq21}) and (\ref{eq32}) that
	
	\begin{align}
		\sum_{n=0}^\infty a_{p-4}(n)q^n&=\dfrac{1}{f_1f_2^{p-5}}\equiv \dfrac{1}{f_{2p}}\cdot\dfrac{f_2^5}{f_1^2}\cdot f_1\nonumber\\
		&\equiv\dfrac{1}{f_{2p}}\sum_{m,\,n=-\infty}^\infty (-1)^{m+n}(3m+1)q^{3m^2+2m+n(3n+1)/2}\pmod{p}.\label{eq33}
	\end{align}
	
	We now consider the equation 
	\[3m^2+2m+\dfrac{n(3n+1)}{2}\equiv l\pmod{p},\]
	which can be written as 
	\begin{align}
		2(6m+2)^2+(6n+1)^2\equiv 3(8l+3)\equiv 0\pmod{p}.\label{eq34}
	\end{align}
	Since $p\equiv 5,7\pmod{8}$, $-2$ is a quadratic nonresidue modulo $p$, so the solution of equation (\ref{eq34}) is $6m+2\equiv 6n+1\equiv 0\pmod{p}$. Thus, we get $3m+1\equiv 0\pmod{p}$, and extracting the terms containing $q^{pn+l}$ on both sides of equation (\ref{eq33}), dividing by $q^l$, and then replacing $q^p$ with $q$ yield equation (\ref{eq31}).
\end{proof}

We next prove the analogous result for primes $p\equiv 3,7\pmod{8}$, which may be seen as a generalization of equation (\ref{eq12}) in Theorem \ref{thm11}. \\

\begin{theorem}\label{thm32}
	Let $p\geq 7$ be a prime with $p\equiv 3,7\pmod{8}$. Then for all $n\geq 0$,
	\begin{equation}\label{eq35}
		a_{p-6}\left(pn+\dfrac{13(p^2-1)}{24}\right)\equiv 0\pmod{p}.
	\end{equation}
\end{theorem}

\begin{proof}
	Since $f_{2p}\equiv f_2^p\pmod{p}$, 
	\begin{align}\label{eq36}
		\sum_{n=0}^\infty a_{p-6}(n)q^n=\dfrac{1}{f_1f_2^{p-7}}\equiv \dfrac{1}{f_{2p}}\cdot\dfrac{f_2^7}{f_1}\pmod{p}.
	\end{align}
	Let 
	\[\sum_{n=0}^\infty A(n)q^n := \dfrac{f_2^7}{f_1}.\]
	Then we have the following identity \cite[p. 223]{ahl}
	\begin{align}\label{eq37}
		A\left(pn+\dfrac{13(p^2-1)}{24}\right)=\epsilon p^2A\left(\dfrac{n}{p}\right),
	\end{align}
	where $p\equiv 7, 11\pmod{12}$ is a prime and 
	\[\epsilon=\begin{cases}
		1 & \text{ if }p\equiv 7\pmod{8},\\
		-1 & \text{ if }p\equiv 3\pmod{8}.
	\end{cases}\]
	As $p\geq 7$ and $p\equiv 3,7\pmod{8}$, we know that $p\equiv 7,11,19,23\pmod{24}$, so $p\equiv 7, 11\pmod{12}$. Thus, applying equation (\ref{eq37}) to equation (\ref{eq36}) yields
	\[a_{p-6}\left(pn+\dfrac{13(p^2-1)}{24}\right)\equiv A\left(pn+\dfrac{13(p^2-1)}{24}\right)\equiv 0\pmod{p}\]
	for any $n\geq 0$, completing the proof of equation (\ref{eq35}).
\end{proof}

\noindent {\bf Acknowledgement.} The author would like to thank the anonymous referee for giving insightful comments that improved the contents of this paper and for bringing the paper \cite{ahl} to his attention, which led to the proof of Theorem \ref{thm32}.

\end{document}